\documentclass[11pt]{article}
\usepackage{graphicx}
\usepackage{mathptmx}
\usepackage{latexsym,amsmath,amssymb,amsfonts,amsthm, amsbsy}

\newcommand{\la}[1]{\mbox{\large $#1$}}
\newcommand{\La}[1]{\mbox{\Large $#1$}}
\newcommand{\LA}[1]{\mbox{\LARGE $#1$}}

\newfont{\lie}{eufm10 at 12pt}
\newfont{\field}{msbm10 at 11pt}

\begin{document}
\title{ Stable $3$-spheres in  $\mathbb{C}^{3}$}
\author{ Isabel M.C.\ Salavessa}
\date{}
\protect\footnotetext{\!\!\!\!\!\!\!\!\!\!\!\!\! {\bf MSC 2000:}
Primary:  53C42; 53C38; 58E35. Secondary:  35J20; 49R50\\
{\bf ~~Key Words:} Stability,  Parallel
Mean curvature,  Calibration,  Cauchy-Riemann inequality,
Spherical harmonics.\\[1mm]
Partially supported by Funda\c{c}\~{a}o para a Ci\^{e}ncia e Tecnologia,  
through
the plurianual project of Centro de F\'{\i}sica das Interac\c{c}\~{o}es Fundamentais and
 programs  PTDC/MAT/101007/2008,  PTDC/MAT/118682/2010.
}
\maketitle ~~~\\[-10mm]
{\footnotesize Centro de F\'{\i}sica das Interac\c{c}\~{o}es
Fundamentais, Instituto Superior T\'{e}cnico, Technical University
of Lisbon, Edif\'{\i}cio Ci\^{e}ncia, Piso 3, Av.\ Rovisco Pais,
1049-001 Lisboa, Portugal;~ isabel.salavessa@ist.utl.pt}\\[5mm]
{\small {\bf Abstract:}
By only using  spectral theory of the Laplace operator on  spheres,
we prove that  the unit 3-dimensional sphere of a 2-dimensional complex 
subspace of $\mathbb{C}^3$  is  an $\Omega$-stable submanifold with parallel
mean curvature, when $\Omega$ is the K\"{a}hler calibration of rank $4$
of $\mathbb{C}^3$.}\\[1mm]

\textbf{1. Introduction}\\
\par
In 2000, Frank Morgan introduced the notion of multi-volume for an
$m$-dimensional submanifold $M$ of a Euclidean space $\mathbb{R}^{m+n}$,
as a volume enclosed by orthogonal projections onto axis $(m+1)$-planes.
He characterized stationary submanifolds for the area functional
with prescribed multi-volume as submanifolds with mean curvature
vector  $H$ prescribed  by a constant multivector 
$\xi\in \wedge_{m+1}\mathbb{R}^{m+n}$, namely $H=\xi\lfloor \vec{S}$, 
where $\vec{S}$ is the unit tangent plane of $M$, and proved  the existence
of a minimizer among rectifiable currents, as well as their regularity under  
general conditions of the boundary. In this setting, a  question  has arisen 
on conditions  for $\|H\|$ to be constant.  
In (Salavessa, 2010) we extended the variational characterization of 
hypersurfaces with constant mean curvature  $\|H\|$ to submanifolds  
with higher codimension, when the ambient space 
is any Riemannian manifold $\bar{M}^{m+n}$, as discovered 
by Barbosa, do Carmo and Eschenburg (1984, 1988) for the case $n=1$.
This generalization amounts on defining an ``enclosed'' $(m+1)$-volume of an
$m$-dimensional immersed submanifold $F:M^m\to \bar{M}^{m+n}$, $m\geq 2$,
 as the $\Omega$-volume defined by
each  one-parameter variation family $F(x,t)=F_t(x)$ of $F(x,0)=F(x)$, 
 where $\Omega$
is a semi-calibration on the ambient space $\bar{M}$, that is, 
  an $(m+1)$-form $\Omega$ which satisfies
$|\Omega(e_0,e_1, \ldots, e_{m})|\leq 1$, for any orthonormal system
$e_i$ of $T\bar{M}$.
A submanifold with calibrated extended tangent space $H\oplus TM$
is a critical point of the functional area, for compactly
supported $\Omega$-volume preserving
variations, if and only if it has  constant mean curvature $\|H\|$.
In this case we have $H=\|H\|\, \Omega\, \lfloor \vec{S}$.
From a deeper inspection of this proof, one can see that the initial assumption
of calibrated extended tangent space can be dropped, since it  will appear as
a consequence of being a critical point itself. 
This will be explained in detail
in a future paper, and also its relations with  Morgan's formalism.
Assuming  that $M$ has parallel mean curvature $H$,  a second variation
 is then computed,
and its non-negativeness defines  stability of $M$. This corresponds
to the non-negativeness of the quadratic form associated with
the $L^2$-self-adjoint  $\Omega$-Jacobi operator
$\mathcal{J}_{\Omega}(W)=\mathcal{J}(W)+m\|H\|C_{\Omega}(W)$, acting 
on  sections in the twisted normal bundle $H^1_{0,T}(NM)=
{\cal F}\oplus H^1_0(E)$, where the set   ${\cal F}$ of $H^1_0$-functions
with zero mean value  is identified with 
the set of  sections of the form $f\nu$, with $f\in {\cal F}$
and $\nu=H/\|H\|$, and where $E$ is the orthogonal complement of $\nu$ in
the normal bundle. 
This Jacobi operator
is the usual one, but   with an extra term, namely a multiple of 
a first order differential operator $C_{\Omega}(W)$ that depends on $\Omega$.
The twisted normal bundle is the $H^1$-completion of the 
vector space generated by
the set $\mathcal{F}_{\Omega}$ of compactly supported
infinitesimal  $\Omega$-volume preserving variations, and, in general, 
 we do not know whether it is larger than $\mathcal{F}_{\Omega}$ itself. 
Thus, $\Omega$-stability implies that the area
functional of $F_t$ decreases when $t$ approaches $t_0=0$, for
 any family of $\Omega$-volume preserving variations $F_t$ of $F$, but we 
do not know whether the converse also holds always.
In case the ambient space is
the Euclidean space $\mathbb{R}^{m+n}$, then a  unit $m$-sphere of an 
$\Omega$-calibrated Euclidean subspace 
$\mathbb{R}^{m+1}$ of $\mathbb{R}^{m+n}$  is $\Omega$-stable if and
only if, for any $(n-1)$-tuple of  functions 
$f_{\alpha}\in C^{\infty}(\mathbb{S}^m)$, $2\leq \alpha \leq n$, the following 
integral inequality holds:
\begin{equation}
\sum_{\alpha<\beta}-2m\int_{\mathbb{S}^m}
f_{\alpha}\xi(W_{\alpha},W_{\beta})(\nabla f_{\beta})dM\leq 
\sum_{\alpha}\int_{\mathbb{S}^m}\|\nabla f_{\alpha}\|^2dM, 
\end{equation}
where  $W_{\alpha}$ is a fixed  global parallel orthonormal (o.n.) frame
of $\mathbb{R}^{n-1}$, the orthogonal complement of $\mathbb{R}^{m+1}$
spanned by $\mathbb{S}^m$, and 
$\xi$
is the $ T^*\mathbb{S}^m$-valued 2-form on $\mathbb{R}^{n-1}_{/\mathbb{S}^{m}}$
$$\xi(W,W')(X)=\Omega(W, W',*X),\quad W,W'\in \mathbb{R}^{n-1}, X\in  T^*\mathbb{S}^m$$
where $*: T\mathbb{S}^m\to\wedge^{m-1}T\mathbb{S}^m$ is the star operator.
If (1)  holds and 
\begin{equation}
\bar{\nabla}_{W}\Omega(W, e_1, \ldots, e_m)=0, 
\quad \forall W\in N\mathbb{S}^m, 
\end{equation} 
where $e_i$ is an o.n.\ frame of $T\mathbb{S}^m$ , 
then in (Salavessa, 2010, proposition 4.5) we have shown that for each 
$\alpha<\beta$,
$\xi(W_{\alpha},W_{\beta})$ must be co-exact as a 1-form on $\mathbb{S}^m$,
that is, 
$$
\xi_{\alpha\beta}:=\xi(W_{\alpha},W_{\beta})=\delta \omega_{\alpha\beta},
$$
for some
globally defined 2-form $ \omega_{\alpha\beta}$ on $\mathbb{S}^m$. 
This is the case
when  $\Omega$ is a parallel $(m+1)$-form
on $\mathbb{R}^{m+n}$.
Using these forms $\omega_{\alpha\beta}$,  the stability condition 
(1) is
translated into the \em long $\Omega$-Cauchy-Riemannian integral 
inequality: \em
\begin{equation}
\sum_{\alpha<\beta}-2m\int_{\mathbb{S}^m}\omega_{\alpha\beta}(\nabla f_{\alpha},\nabla f_{\beta})
dM\leq \sum_{\alpha}\int_{\mathbb{S}^m}\|\nabla f_{\alpha}\|^2dM.
\end{equation}
If we fix $\alpha<\beta$, and set $f=f_{\alpha}$,
$h=f_{\beta}$, and $f_{\gamma}=0$ $\forall \gamma
\neq \alpha,\beta$, (1) reduces to
\begin{equation}
-2m\int_{\mathbb{S}^m}f\xi_{\alpha\beta}(\nabla h)
dM\leq \int_{\mathbb{S}^m}\|\nabla f\|^2dM +
 \int_{\mathbb{S}^m}\|\nabla h\|^2dM, 
\end{equation}
and if we  replace $f$ by $cf$, and $h$ by $c^{-1}h$, where $c^2=
\|\nabla h\|_{L^2}/\|\nabla f\|_{L^2}$, 
then we obtain the  corresponding  equivalent 
\em short $\Omega$-Cauchy-Riemannian,
integral inequality \em
\begin{eqnarray}
-m\int_{\mathbb{S}^m}\omega_{\alpha\beta}(\nabla f, \nabla h)
dM\leq \sqrt{\int_{\mathbb{S}^m}\|\nabla f\|^2dM} 
\sqrt{ \int_{\mathbb{S}^m}\|\nabla h\|^2dM}, 
\end{eqnarray}
holding for all functions $f,h\in C^{\infty}(\mathbb{S}^m)$.

The $\Omega$-stability  of a submanifold with calibrated extended tangent 
space and parallel mean curvature
depends on the curvature of the ambient space and on the
calibration $\Omega$ (Salavessa, 2010). It always holds on Euclidean spheres 
if $C_{\Omega}$ vanish. 
This last condition is equivalent to the condition (2) and 
$\xi\equiv 0$ ((Salavessa, 2010), Lemma 4.4). In the case $n=2$ the later condition is satisfied,
but for $n \geq 3$ the operator $C_{\Omega}$ 
may not vanish for spheres, even if $\Omega$ is parallel. If
$C_{\Omega}$ does not vanish, spheres of calibrated vector subspaces 
 may not be $\Omega$-stable. 

We first consider  $\Omega$  any parallel $(m+1)$-form
on $\mathbb{R}^{m+n}$.
Laplace spherical harmonics of $\mathbb{S}^m$ of degree $l$ are
the eigenfunctions for the closed  eigenvalue problem with respect to
the  Laplacian operator
corresponding to the eigenvalue $\lambda_l=l(l+m-1)$, and
they are just the harmonic homogeneous polynomial functions 
of degree $l$ of $\mathbb{R}^{m+1}$ restricted to $\mathbb{S}^m$.
We denote by $E_{\lambda_l}$ the finite-dimensional subspace of
$H^1( \mathbb{S}^m)$ spanned by these $\lambda_l$-eigenfunctions.
In the first theorem we show how each 1-form
$\xi_{\alpha\beta}$ transforms a spherical harmonic $f$ into another
spherical harmonic $h$:\\[1mm]
{\bf Theorem 1.1.}  \em If $\Omega$ is parallel,  then for each
$f\in E_{\lambda_l}$,  $h=\xi_{\alpha\beta}(\nabla f)$ is also in
$E_{\lambda_l}$, and it is $L^2$-orthogonal to $f$. \em \\

In this paper we study the stability of
the unit 3-sphere of a 2-dimensional complex subspace
of $\mathbb{C}^3$ with respect to the K\"{a}hler calibration. In this case
$C_{\Omega}$ does not vanish.
Let $\varpi$ be the K\"{a}hler form of $\mathbb{C}^3=\mathbb{R}^6$,
and $\Omega$ the  K\"{a}hler calibration of rank 4,
$$ \varpi=dx^{12}+ dx^{34}+dx^{56}, \quad~~~~~
\Omega=\frac{1}{2}\varpi^2.$$
The unit sphere of $\mathbb{R}^4\times\{0\}$
is immersed into $\mathbb{R}^6=\mathbb{C}^3$, by the inclusion map
$\phi=(\phi_1,\ldots,\phi_4,0):\mathbb{S}^3 \to \mathbb{C}^3$.
We have only one of those 1-forms
$$\xi:=\xi_{56}= *(d\phi^1\wedge d\phi^2+ d\phi^3\wedge d\phi^4)=
\phi^1 d\phi^2-\phi^2 d\phi^1+  \phi^3 d\phi^4-\phi^4 d\phi^3,$$
and  $\xi= \delta \omega,$ with
$\omega=\frac{1}{2}*\xi=\frac{1}{2}
(d\phi^1\wedge d\phi^2+ d\phi^3\wedge d\phi^4)=\frac{1}{2}\phi^*\varpi$.
Our main theorem is the following:\\[2mm]
{\bf Theorem 1.2.}  \em  Three-dimensional spheres of $\mathbb{C}^2$
are $\Omega$-stable submanifolds  of $\mathbb{C}^3$ 
with parallel mean curvature,  where
$\Omega= \frac{1}{2}\varpi^2$ is
the K\"{a}hler calibration of rank $4$. \em \\[2mm]
The Cauchy-Riemann inequality version of the $\Omega$-stability is described in
the corollary:\\[2mm]
{\bf Corollary 1.1.} \em 
The  Cauchy-Riemann inequality
$$-\int_{\mathbb{S}^3} \varpi(\nabla f,\nabla h)dM
\leq \frac{2}{3}\sqrt{\int_{\mathbb{S}^3}\|\nabla f\|^2dM}
 \sqrt{\int_{\mathbb{S}^3}\|\nabla h\|^2dM}$$
holds for any smooth functions $f$ and $h$ of
$~\mathbb{S}^3$,
with equality if and only if $f, h\in E_{\lambda_1}$, with
$f=\sum_i\mu_i\phi_i$ and $h=\sum_i\sigma_i\phi_i$, where
$\sigma_2=-\mu_1$, $ \sigma_1=\mu_2$, $\sigma_4=-\mu_3$, $\sigma_3=\mu_4$.\em 
\\[1mm]

Finally, we state that  the 3-sphere is the unique smooth closed submanifold 
that solves the $\Omega$-isoperimetric problem among a certain class of
immersed submanifolds:\\[2mm]
{\bf Theorem 1.3.} \em
The unit 3-sphere of a complex 2-dimensional subspace
of $\mathbb{C}^3$ is the unique closed immersed 3-dimensional submanifold
$\phi: M\to \mathbb{C}^{3}$
with parallel mean curvature, trivial normal bundle,
 and complex extended  tangent
space $H\oplus TM$, that is $\Omega$-stable for the K\"{a}hler
calibration of rank 4, and satisfies  the inequality 
$$\int_M S(2 + h \|H\|)dM \leq 0,$$ 
where $h$ and $S$ are the height functions
$ h=\langle \phi, \nu\rangle$ and
$S= \sum_{ij} \langle \phi, (B(e_i,e_j))^F \rangle B^{\nu}(e_i,e_j).$
\em \\[2mm]
{\bf Remark. } On a closed K\"{a}hler manifold $(M,J)$ with K\"{a}hler form
$\varpi(X,Y)=g(JX,Y)$, if $f,h:M\to \mathbb{R}$ are smooth functions, 
then by the Cauchy-Schwarz inequality,
$$\left|\int_M\varpi(\nabla f, \nabla h)dM\right| \leq
\sqrt{\int_{M}\|\nabla f\|^2dM}
 \sqrt{\int_{M}\|\nabla h\|^2dM}, $$
with equality if and only if $\nabla h=\pm J\nabla f$, or equivalently
$f\pm ih:M\to \mathbb{C}$ is a holomorphic function. If this is the case,
then $f$ and $h$ are constant functions. On the other hand, globally
defined functions, 
sufficiently close to holomorphic functions defined on a sufficiently large
open set, are expected to satisfy an almost equality.
This is not the case of $\mathbb{S}^3$, which is
not a complex manifold, and  somehow explains the  coefficient $2/3$ 
 in  Corollary 1.1.\\[2mm]
{\bf  Remark.}  In the case of $3$-spheres in $\mathbb{C}^3$ we have only 
one form $\xi_{\alpha\beta}$, that is, the long Cauchy-Riemann inequality
is the short one.  We wonder if a general proof of 
short Cauchy-Riemann inequalities can be allways obtained
  for Euclidean $m$-spheres on $\mathbb{R}^{m+n}$, by  using the
 spectral theory of spheres, when
$\Omega$ is any parallel calibration.
 Note that (4) is immediately satisfied for
$f,h\in E_{\lambda_l}$, if $\lambda_l\geq m^2$, that is $l\geq m$, so
it remains to consider the cases $l\leq m-1$.
For 3-spheres we have to 
consider  polynomial functions up to order $l=2$, while for 2-spheres
we have  to consider only the case $l=1$.
A related remark is given in 
the end of section 3.\\[1mm]

\textbf{2. Preliminaries}\\

We consider an oriented  Riemannian manifold $M$ of dimension $m$, 
with Levi-Civita connection $\nabla$ and Ricci tensor $Ricci^M:TM \to TM$.
In what follows $e_1, \ldots, e_m$
denotes a local direct o.n. frame.\\[2mm]
{\bf Lemma 2.1.} \em 
Let $\xi$ be a co-exact 1-form on a Riemannian manifold $M$,
with $\xi=\delta \omega$, where $\omega$ is a 2-form. Then for any
function $f\in C^{2}(M)$, 
$$\xi(\nabla f)=div(\nabla^{\omega}f),$$
where $\nabla^{\omega}f= \sum_i\omega(\nabla f,e_i)e_i~$.
 Moreover, for any $f,h\in C^{\infty}_0(M)$
$$\int_M f \xi(\nabla h)dM =\int_{M}\omega(\nabla f, \nabla h)dM=-
\int_M h \xi(\nabla f)dM.$$ \em \\[2mm]
{\bf Proof.}  We may assume at a point $x_0$, $\nabla e_i=0$. Then at $x_0$
\begin{eqnarray*}
\xi(\nabla f)&=&\delta \omega(\nabla f)=-\sum_i\nabla_{e_i}\omega(e_i, 
\nabla f)
=\sum_i-\nabla_{e_i}(\omega(e_i, \nabla f))+ \omega(e_i, \nabla_{e_i}\nabla f)\\
&=& div (\nabla^{\omega} f)+ \sum_{ij}Hess f(e_i,e_j)\omega(e_i,e_j).
\end{eqnarray*}
The  last equality proves the first equality of the lemma, because
$Hess f(e_i,e_j)$ is symmetric on $i,j$ and $\omega(e_i,e_j)$ is skew-symmetric.
The other equalities of the lemma follow from  $div(f X)= 
\langle \nabla f, X\rangle+f div(X)$, holding for any vector field $X$
and function $f$.
\qed\\

The $\delta$ and  star operators  acting on 
$p$-forms on an oriented  Riemannian $m$-manifold $M$ satisfy
$\delta=(-1)^{mp+m+1}*d*$, $**=(-1)^{p(m-p)}Id$, and for a 
1-form $\xi$  the DeRham Laplacian $\Delta$ 
and the rough Laplacian $\bar{\Delta}$ are related
by the following formulas 
$$
 \begin{array}{l}
\Delta \xi (X)= (d\delta + \delta d)\xi (X)= 
-\bar{\Delta}\xi (X) +\xi(Ricci^M(X)),\\[1mm]
\bar{\Delta}\xi (X)= trace \nabla^2 \xi(X)= 
\sum_i \nabla_{e_i}\nabla_{e_i}\xi(X)
-\nabla_{\nabla_{e_i}e_i}\xi(X). \end{array}
$$
If $\xi=\delta \omega$, then  $\delta \xi=0$, and so
$\Delta \xi(X)=\delta d \xi (X)= -\sum_i \nabla_{e_i}(d\xi) (e_i,X)$.
We also recall the following well-known formula (see e.g. Salavessa \& Pereira do Vale 
(2006)) for $f\in C^{\infty}(M),$
$$(\bar{\Delta}df)(X)=\sum_i\nabla^2_{e_i,e_i}df (X)
=g(\nabla (\Delta f), X) +df(Ricci^M(X)).$$
Thus, 
\begin{equation}
\begin{array}{l}
\bar{\Delta} (\nabla f) = \nabla (\Delta f)+ Ricci^M(\nabla f),\\[1mm]
(\bar{\Delta }\xi)(\nabla f)=-(\delta d\xi)(\nabla f) +\xi(Ricci^M(\nabla f)).
\end{array}
\end{equation}

Now we suppose that $M$ is an immersed oriented
 hypersurface of a Riemannian manifold 
${M'}$, with Riemannian metric $\langle, \rangle$, 
defined by an immersion $\phi:M\to {M'}$ with unit normal
$\nu$, second fundamental form $B$ and corresponding Weingarten operator $A$
in the $\nu$ direction, given by
$$B(e_i,e_j)= \langle A(e_i),e_j\rangle 
=\langle {\nabla'}_{e_i}e_j,\nu\rangle =
-\langle e_j, \nabla'_{e_i}\nu\rangle, $$
where ${\nabla'}$ denotes the Levi-Civita connection on ${M'}$. 
The scalar mean curvature of $M$ is given by
$$ H=\frac{1}{m} Trace\,  B= \sum_i\frac{1}{m} B(e_i,e_i).$$
The curvature operator of ${M'}$,  
${R'}(X,Y,Z,W)=\langle
-{\nabla'}_X{\nabla'}_YZ+{\nabla'}_Y{\nabla'}_XZ
+{\nabla'}_{[X,Y]}Z, W\rangle$,
can be seen as a self-adjoint operator of wedge bundles
 ${R'}: \wedge^2 T{M'} \to \wedge^2 T{M'}$,
$$\langle {R'}(u\wedge v), z\wedge w\rangle
={R'}(u,v,z,w),$$
and so $R'(u\wedge v)=\sum_{i<j}  R'(u, v, e_i, e_j)e_i\wedge e_j$, 
where 
$$<u\wedge v, z\wedge w>= det\left[\begin{array}{cc}
\langle u, z\rangle & \langle u,w\rangle\\
\langle v, z\rangle & \langle v, w\rangle \end{array}\right].$$
In what follows, we suppose that
 $\hat{\xi}$ is  a parallel $(m-1)$-form on ${M'}$, and $\xi$ is given by
$$\xi=*\phi^*\hat{\xi}$$
where $*$ is the star operator on $M$. In this case
$\xi$ is obviously co-closed, but
not necessarily co-exact.
We employ the usual inner products in $p$-forms and morphisms.\\[2mm]
{\bf Lemma 2.2.} \em  Assume $m\geq 3$. Then for all $i,j$ 
$$\begin{array}{l}
(\nabla_{e_i}\xi) (e_j)=
\sum_{k}-B(e_i,e_k)\hat{\xi}(\nu, *(e_k\wedge e_j))=
-\hat{\xi}(\nu, *(A(e_i)\wedge e_j)),\\[1mm]
\Delta\xi(e_j)= \delta\,  d \xi(e_j) =\hat{\xi}\left( \nu, *\la{(}e_j\wedge (
 m\nabla H -[Ricci^{{M'}}(\nu)]^T)\la{)} + {R'}(e_j\wedge \nu)\right) 
+ \xi(\Theta_B(e_j)),
\end{array}$$
where   $[Ricci^{{M'}}(\nu)]^T=\sum_k Ricci^{{M'}}(\nu, e_k)e_k$
and  $\Theta_B:T{M}\to T{M}$ is the morphism given  by,
$\Theta_B=\|B\|^2 Id + mH A -2A^2$. \em \\[2mm]
{\bf  Proof.}  We fix a point $x_0\in M$
and take $e_i$ a local o.n. frame s.t.
$\nabla e_i(x_0)=0$. 
  We will compute $d\xi(e_i,e_j)$,
at $x$ on a neigbourhood of $x_0$.
Recall that
for any $p$-form $\sigma$, we have  $*\sigma=\sigma *$,  where the star operator on the
r.h.s. can be seen as acting on $\wedge^{m-p}TM$, with
$*e_i=(-1)^{i-1}e_1\wedge \ldots \wedge \hat{e}_i\wedge \ldots e_m$,  and 
for $i<j$, $*(e_i\wedge e_j)=(-1)^{i+j-1}e_1\wedge\ldots \wedge \hat{e}_i
\wedge \ldots\wedge  \hat{e}_j \wedge \ldots \wedge e_m$.
Using the fact that $\hat{\xi}$ is a parallel form on $M'$, 
we have  for $x$ near $x_0$,
\begin{eqnarray*}
\begin{array}{rl}
\nabla_{e_i}(\xi (e_j))=&\sum_{k\neq j}(-1)^{j-1}\hat{\xi}(e_1,
\ldots ,{\nabla'}_{e_i}e_k, \ldots, \hat{e}_j, \ldots, e_m)\\[1mm]
=& \sum_{k<j}(-1)^{k+j}\hat{\xi}( {\nabla'}_{e_i}e_k, e_1,
\ldots ,\hat{e}_k, \ldots, \hat{e}_j, \ldots, e_m)\\
&+\sum_{k>j}(-1)^{k+j-1}\hat{\xi}( {\nabla'}_{e_i}e_k, e_1,
\ldots ,\hat{e}_j, \ldots, \hat{e}_k, \ldots, e_m)\\[1mm]
=&\sum_{k<j}-\langle \nabla_{e_i}e_k,e_j\rangle \hat{\xi}(*e_k)
-B(e_i,e_k)\hat{\xi}(\nu, *(e_k\wedge e_j))\\
&+\sum_{k>j}-\langle \nabla_{e_i}e_k,e_j\rangle \hat{\xi}(*e_k)
+B(e_i,e_k)\hat{\xi}(\nu, *(e_j\wedge e_k))\\[1mm]
=& \xi(\nabla_{e_i}e_j)+  
\sum_{k\neq j}-B(e_i,e_k)\hat{\xi}(\nu, *(e_k\wedge e_j)).
\end{array}
\end{eqnarray*}
Hence, 
$(\nabla_{e_i}\xi) (e_j) =
\sum_{k\neq j}-B(e_i,e_k)\hat{\xi}(\nu, *(e_k\wedge e_j))$, which
 proves the first sequence of equalities of the lemma.
Now,
\begin{eqnarray*}
 d\xi(e_i,e_j) &=& (\nabla_{e_i}\xi) (e_j)-(\nabla_{e_j}\xi) (e_i)\\
&=&
\sum_{k\neq j}-B(e_i,e_k)\hat{\xi}(\nu, *(e_k\wedge e_j))
+\sum_{k\neq i}B(e_j,e_k)\hat{\xi}(\nu, *(e_k\wedge e_i)),
\end{eqnarray*}
and by Codazzi's equation,
$$ \begin{array}{l}
(\nabla_{e_i}B)(e_j,e_k)= (\nabla_{e_j}B)(e_i,e_k)-{R'}(e_i,e_j,e_k,\nu)\\
\sum_i (\nabla_{e_i}B)(e_i,e_k)=m\nabla_{e_k}H- Ricci^{{M'}}(e_k,\nu).
\end{array}$$
Note that $B_{ik}=(\nabla_{e_j}B)(e_i,e_k)$ is a symmetric
matrix, and if we define $A_{ki}=\hat{\xi}(\nu,* (e_k\wedge e_i))$
(valuing zero if $k=i$), then $A_{ik}$ is skew-symmetric.
Thus,  $\sum_{k\neq i} B_{ik}A_{ki}=\sum_{k,i}B_{ik}A_{ki}=0.$
Furthermore, if we set
$C_{ik}=-{R'}(e_i,e_j,e_k,\nu)$, then
$C_{ik}-C_{ki}={R'}(e_k,e_i,e_j,\nu)$. Hence, 
$$\sum_{i}\sum_{k\neq i} C_{ik}A_{ki}=\sum_{ik} C_{ik}A_{ki}=
\sum_{ik} \frac{1}{2}((C_{ik}+ C_{ki})+ (C_{ik}- C_{ki}))A_{ki}
=\sum_{ki}\frac{1}{2}{R'}(e_k,e_i,e_j,\nu)A_{ki}.$$
Therefore,  for each $j$, at $x_0$ 
\begin{eqnarray*}
\lefteqn{-\delta d \xi (e_j)
= \sum_i \nabla_{e_i}(d\xi(e_i,e_j))}\\
&=& \sum_{k\neq j}\sum_i-(\nabla_{e_i}B)(e_i,e_k)\hat{\xi}(\nu,*(e_k\wedge e_j))
-B(e_i,e_k)\nabla_{e_i}(\hat{\xi}(\nu,*(e_k\wedge e_j)))\\
&&+ \sum_{k\neq i}\sum_j(\nabla_{e_i}B)(e_j,e_k)\hat{\xi}(\nu,*e_k\wedge e_i))
+B(e_j,e_k)\nabla_{e_i}(\hat{\xi}(\nu, *(e_k\wedge e_i))\\
&=&\sum_{k\neq j} (-m\nabla_{e_k}H + Ricci^{{M'}}(e_k, \nu))
\hat{\xi}(\nu, *(e_k\wedge e_j))
 +\sum_{k,i} \frac{1}{2}
{R'}(e_k,e_i,e_j,\nu)\hat{\xi}(\nu,*(e_k\wedge e_i)) + S
\end{eqnarray*}
where
\begin{eqnarray*}
 \begin{array}{lcl} 
S&=&\sum_i\sum_{k<j}(-1)^{k+j}B(e_i,e_k)\hat{\xi}({\nabla'}_{e_i}\nu,
e_1, \ldots, \hat{e}_k, \ldots, \hat{e}_j, \ldots, e_m)\\
&&+ \sum_i\sum_{k>j}(-1)^{k+j-1}B(e_i,e_k)\hat{\xi}({\nabla'}_{e_i}\nu,
e_1, \ldots, \hat{e}_j, \ldots, \hat{e}_k, \ldots, e_m)\\
&&+\sum_i\sum_{k<i}(-1)^{k+i-1}B(e_j,e_k)\hat{\xi}({\nabla'}_{e_i}\nu,
e_1, \ldots, \hat{e}_k, \ldots, \hat{e}_i, \ldots, e_m)\\
&&+ \sum_i\sum_{k>i}(-1)^{k+i}B(e_j,e_k)\hat{\xi}({\nabla'}_{e_i}\nu,
e_1, \ldots, \hat{e}_i, \ldots, \hat{e}_k, \ldots, e_m)
\end{array}\\
 \begin{array}{lcl}
=&&\sum_{i}\sum_{k<j}-B(e_i,e_k)B(e_i, e_k)\xi(e_j)+
B(e_i,e_j)B(e_i, e_k)\xi(e_k)\\
&&+\sum_{i}\sum_{k>j}B(e_i,e_j)B(e_i, e_k)\xi(e_k)-
B(e_i,e_k)B(e_i, e_k)\xi(e_j)\\
&&+ \sum_{i}\sum_{k<i}B(e_i,e_k)B(e_j, e_k)\xi(e_i)
-B(e_i,e_i)B(e_j, e_k)\xi(e_k)\\
&&+\sum_{i}\sum_{k>i}-B(e_i,e_i)B(e_j, e_k)\xi(e_k)+
B(e_i,e_k)B(e_j, e_k)\xi(e_i).
\end{array}
\end{eqnarray*}
At this point we may assume that at $x_0$ the basis $e_i$ diagonalizes
the second fundamental form, that is, $B(e_i,e_j)=\lambda_{i}\delta_{ij}$.
Then,
\begin{eqnarray*}
\begin{array}{lcl}
 S &=&
\sum_i\sum_{k<j}-\delta_{ik}\lambda_i^2\xi(e_j)
+ \delta_{ij}\delta_{ik}\lambda_i^2\xi(e_k)+\sum_i\sum_{k>j}
\delta_{ij}\delta_{ik}\lambda_i^2\xi(e_k)- \delta_{ik}\lambda_i^2\xi(e_j)\\
&&+\sum_i\sum_{k<i}\delta_{ik}\delta_{jk}\lambda_k^2\xi(e_i)
-\delta_{ii}\delta_{jk}\lambda_i\lambda_j\xi(e_k)
+\sum_i\sum_{k>i}
-\delta_{ii}\delta_{jk}\lambda_i\lambda_j\xi(e_k)
+\delta_{ik}\delta_{jk}\lambda_k^2\xi(e_i)\\[1mm]
&=& \sum_{i<j}-\lambda_i^2\xi(e_j)+ \sum_{i>j}-\lambda_i^2\xi(e_j)
+\sum_{j<i}-\lambda_i\lambda_j\xi(e_j)+ \sum_{j>i}-\lambda_i\lambda_j\xi(e_j)
\\[1mm]
&=& \sum_{i\neq j}-\lambda_i^2\xi(e_j)-\lambda_i\lambda_j\xi(e_j)
=  \sum_{i}-\lambda_i^2\xi(e_j)-\lambda_i\lambda_j\xi(e_j) + (\lambda_j^2
+\lambda_j^2)\xi(e_j)\\[1mm]
&=& -\|B\|^2\xi(e_j)-mH\xi(A(e_j))+2\xi(A^2(e_j)),
\end{array}
\end{eqnarray*}
and the second sequence of equalities of the lemma is proved.\qed\\[2mm]
If we suppose that $\Theta_B= \mu(x)Id$, taking 
$e_i$ a diagonalizing o.n. basis of the second fundamental form,
 $B(e_i,e_j)=\lambda_i \delta_{ij}$, then
 each $\lambda_i$ satisfies the quadratic equation 
$$2\lambda_i^2-mH\lambda_i + (\mu-\|B\|^2)=0,$$
which implies that we have at most two distinct possible principal curvatures
$\lambda_{\pm}$.
Moreover, from the above equation, summing over $i$, we derive that $\mu(x)$ must satisfy $\mu(x)=\frac{m-2}{m}\|B\|^2+ m H^2$, and so
$$\lambda_{\pm}=\frac{1}{4}\LA{(}mH\pm \sqrt{\frac{16}{m}\|B\|^2+ m(m-8)H^2}
\LA{)}.$$
Note that, from $\|B\|^2\geq m \|H\|^2$, we have
 $\frac{16}{m}\|B\|^2+ m(m-8)H^2\geq (m-4)^2H^2$, and so
there are one or two distinct principal curvatures. If
 $M$ is totally umbilical, then $\|B\|^2=mH^2$ and $\mu=2(m-1)\|H\|^2$.
 The previous lemma leads to the following conclusion:\\[1mm]
{\bf Lemma 2.3.} \em 
Assuming ${M'}=\mathbb{R}^{m+1}$, $m\geq 3$, and taking $M$  a hypersurface
with constant mean curvature,  with  $\Theta_B= \mu(x)Id$, where
$\mu(x)$ is a smooth function on $M$,  we get
$\mu(x)=\frac{m-2}{m}\|B\|^2+ m H^2$ and 
$$ \Delta \xi =\mu \xi.$$
Furthermore,
$\xi$ is an eigenform for the DeRham Laplacian operator, that is $\mu(x)$
is constant,  if and only if
$\|B\|$ is constant.

In case $M$ is a unit $m$-sphere $ \mathbb{S}^m$,  then $\Theta_B=\mu Id$,
with $\mu=2(m-1)$, and taking
$\nu_x=-x$ as unit normal, then, at each $x\in \mathbb{S}^m$,
$$\begin{array}{l}
(\nabla_{e_i}\xi) (e_j) = \hat{\xi}(x, *(e_i\wedge e_j))\\
d\xi(e_i,e_j) = 2\hat{\xi}(x, *(e_i\wedge e_j))\\
\Delta\xi=\delta d\xi =2(m-1)\xi.\end{array}$$ \em 
{\bf Lemma 2.4.}  \em If $f\in C^{\infty}(\mathbb{S}^{m})$, then
$\Delta(\xi(\nabla f))=\xi(\nabla \Delta f).$ \em \\[2mm]
\noindent
{\bf  Proof. }  We fix a point $x_0\in \mathbb{S}^m$
and take $e_i$ a local o.n. frame of the sphere s.t.
$\nabla e_i(x_0)=0$. Let $f\in C^{\infty}(\mathbb{S}^{m})$.
The following computations are at $x_0$. Using the above formulas (6)
and previous lemma,  we have
\begin{eqnarray*}
\Delta(\xi(\nabla f))&=& \sum_i{\nabla_{e_i}}(\nabla_{e_i}(\xi(\nabla f)))=
\sum_i{\nabla_{e_i}}\La{(}(\nabla_{e_i}\xi)(\nabla f) 
+ \xi(\nabla_{e_i}\nabla f)\La{)}\\
&=& (\bar{\Delta}\xi )(\nabla f) + 2(\nabla_{e_i}\xi) (\nabla_{e_i}\nabla f)
+\xi(\nabla_{e_i}\nabla_{e_i}\nabla f)\\
&=& -2(m-1)\xi(\nabla f) +\xi(\nabla \Delta f) + 2(m-1)\xi(\nabla f)
+ \sum_i 2(\nabla_{e_i}\xi)(\nabla_{e_i}\nabla f).
\end{eqnarray*}
Since  $Hess\, f(e_i,e_j) $ is symmetric in $ij$ and
by  Lemma  2.3, $(\nabla_{e_i}\xi)(e_j)$ is skew-symmetric, we have 
$$ \sum_i (\nabla_{e_i}\xi)(\nabla_{e_i}\nabla f)=\sum_{ij}
Hess\, f(e_i,e_j) (\nabla_{e_i}\xi)(e_j)=0,$$
and the lemma is proved. \qed \\

\textbf{3. Proof of  Theorem 1.1.}\\

We denote by $\nabla$ the Levi-Civita connection of $\mathbb{S}^m$
induced by the flat connection  $\bar{\nabla}$ of $\mathbb{R}^{m+n}$.
We are considering a parallel calibration $\Omega$ on  $\mathbb{R}^{m+n}$.
We fix $\alpha <\beta$ and  define the 1-form on $\mathbb{S}^{m}$
$$\xi=\xi(W_{\alpha},W_{\beta})= * \phi^*\hat{\xi}=\delta\omega,$$
where $\hat{\xi}=\hat{\xi}_{\alpha\beta}$ and $\omega= \omega_{\alpha\beta}$.

We recall that the eigenvalues of $\mathbb{S}^m$ for the closed Dirichlet 
problem are given by $\lambda_l={l(l+m-1)}$, 
with $l=0,1,2 \ldots$. 
We denote by $E_{\lambda_l}$ the eigenspace of dimension $m_l$
corresponding to the 
eigenvalue  $\lambda_l$,
and by $E_{\lambda_l}^+$ the $L^2$-orthogonal complement of 
the sum of the eigenspaces
$E_{\lambda_i}$, $i=1,\ldots,l-1$, and so it is
the sum of all eigenspaces $E_{\lambda}$ with $\lambda\geq \lambda_l$.
If $f\in E_{\lambda_l}$, and 
$h\in E_{\lambda_s}$, then 
$$\int_{\mathbb{S}^m}fh\, dM=0 ~~\mbox{if}~
l\neq s \mbox{~~~~and~~~~} 
\int_{\mathbb{S}^m}\langle \nabla f, \nabla h\rangle\,  dM=  \delta_{ls}
\lambda_l\int_{\mathbb{S}^m}fh\, dM.$$
There exists an  $L^2$-orthonormal basis $\psi_{l,\sigma}$ 
of $L^2(\mathbb{S}^m)$ 
of eigenfunctions
($1\leq \sigma\leq m_l$).
The Rayleigh characterization of $\lambda_l$ is given by
$$\lambda_l=\inf_{f\in E_{\lambda_l}^{+}}\frac{\int_{\mathbb{S}^m}
\|\nabla f\|^2dM}{\int_{\mathbb{S}^m}f^2dM},$$ 
and the infimum is attained for 
$f\in E_{\lambda_l}$. 
Each eigenspace $E_{\lambda_l}$ is exactly composed
by the restriction to
$\mathbb{S}^m$ of the harmonic 
homogeneous polynomial functions of degree $l$ of $\mathbb{R}^{m+1}$, and
it  has dimension $m_l=\binom{m+l}{m}-\binom{m+l-2}{m}$.
Thus, each eigenfunction
$\psi\in E_{\lambda_l}$ is of the form
$\psi=\sum_{|a|=l} \mu_{a }
\phi^{a}$, where $\mu_{a}$ are
some scalars and
$a=(a_1,\ldots, a_{m+1})$ denotes a multi-index
of length
$|a|=a_1+ \ldots +a_{m+1}=l$ and
$$\phi^{a}=
\phi_1^{a_1}\cdot\ldots\cdot\phi_{m+1}^{a_{m+1}}.$$
From $\nabla \phi_i= \epsilon_i^{\top}$ and  $\sum_i \phi_i^2=1$,
we see that
\begin{equation}
\left\{ \begin{array}{lcl}
\langle \nabla \phi_i,\nabla\phi_j\rangle =\delta_{ij}-
\phi_i\phi_j&&
\|\nabla \phi_i\|^2= 1-\phi_i^2\\[1mm]
\int_{\mathbb{S}^m}\phi^2_idM=\frac{1}{m+1}|\mathbb{S}^m| &&
\int_{\mathbb{S}^m}\|\nabla \phi_i\|^2dM
=\lambda_1\int_{\mathbb{S}^2} \phi_i^2dM =\frac{m}{m+1}|\mathbb{S}^m|.
\end{array}\right.
\end{equation}
We also denote by $\int_{\mathbb{S}^m}\phi^2dM$
any of the integrals $\int_{\mathbb{S}^m}\phi_i^2dM$, $i=1,\ldots, m+1$.
We recall the following:\\[2mm]
{\bf Lemma 3.1.}  \em 
If $P:\mathbb{S}^m\to \mathbb{R}$ is a homogeneous polynomial 
function of degree $l$, then
$$\int_{\mathbb{S}^m}P(x)dM=\frac{1}{\lambda_l}\int_{\mathbb{S}^m}
\Delta^0P(x)dM.$$
In particular,
$$\int_{\mathbb{S}^m}\phi^a dM=\sum_{1\leq i\leq m+1}\frac{a_i(a_i-1)}{l(l+m-1)}
\int_{\mathbb{S}^m}\phi^{a-2\epsilon_i}dM,$$
where the terms $a_i<2$ are considered to vanish.
Thus, if some $a_i$ is odd this integral vanishes. \em \\[1mm]
\noindent 
{\bf  Proof of Theorem 1.1}.
By Lemma 2.4, if $f\in E_{\lambda_k}$ then 
$ \xi(\nabla f)\in E_{\lambda_k}$.
From
$$\int_{\mathbb{S}^m}f \xi(\nabla f)dM=\int_{\mathbb{S}^m}\omega(
\nabla f, \nabla f)dM=0$$
we conclude that $f$ and $h= \xi(\nabla f)$ are $L^2$-orthogonal.\qed\\[4mm]
{\bf  Remark. }   Let us consider  $f,h\in E_{\lambda_l}$, and take the
globally defined vector field of $\mathbb{S}^m$,
$\xi^{\sharp}=\sum_j\xi(e_j)e_j$. From  Lemma 2.2,  we have
$$
\langle \nabla h, \nabla(\xi(\nabla f))\rangle
=-\hat{\xi}(
\nu,*(\nabla h\wedge \nabla f))+ Hess f(\nabla h,\xi^{\sharp}).
$$ 
By Theorem 1.1, 
$\xi(\nabla f)\in  E_{\lambda_l}$ as well. 
The term $Hess f(\nabla h, \xi^{\sharp})$ is a sum of  polynomial functions
of degree $2l-3 + k_{\xi}$ where $k_{\xi}$ depends on $\xi^{\sharp}$,
when expressed in terms of $\phi^i$. 
Let us suppose that all $k_{\xi}$ are even. Then
by Lemma 3.1,
$\int_{\mathbb{S}^m} Hessf(\nabla h, \xi^{\sharp})dM=0$.
Since $\lambda_l\geq m$,
and taking into consideration that $\Omega$ is
a semi-calibration, 
\begin{eqnarray*}
-\int_{\mathbb{S}^m} h\xi(\nabla f)dM &=& -\frac{1}{\lambda_l}
\int_{\mathbb{S}^m} \langle \nabla h,\nabla (\xi(\nabla f))\rangle dM
=
\frac{1}{\lambda_l}
\int_{\mathbb{S}^m}\hat{\xi}(\nu,*(\nabla h\wedge \nabla f))dM\\
&\leq & \frac{1}{\lambda_l}\int_{\mathbb{S}^m}
\|\nabla h\|\, \| \nabla f\|dM\leq \frac{1}{m}
\|\nabla f\|_{L^2}\|\nabla h\|_{L^2}.
\end{eqnarray*}
Thus, in this case the short Cauchy-Riemann inequality holds.
Inspection of $\xi$ must be required for each case of $\Omega$.
A general proof of the short Cauchy-Riemann integral inequality, 
under appropriate conditions 
on $\Omega$, will be developed in a future paper.\\[1mm]

\textbf{4. 3-spheres of $\mathbb{C}^2$ in $\mathbb{C}^3$}\\
\par

In this section we specialize  the Cauchy-Riemann inequalities
for the case $m=n=3$ and for $\mathbb{R}^6= \mathbb{C}^3$
we will consider  the K\"{a}hler calibration 
$\frac{1}{2}\mathbf{\varpi}^2$
that calibrates the complex two-dimensional subspaces, that is,  
$$\Omega= dx^{1234}+ dx^{1256}+ dx^{3456}.$$
Thus, fixing $W_5=\epsilon_5$ and $W_6=\epsilon_6$ we have
$\hat{\xi}:= \hat{\xi}_{56}= dx^{12}+ dx^{34}$, and
$$\xi:=\xi_{56}= *\phi^*\hat{\xi}=*(d\phi^{12}+d\phi^{34}).$$
The volume element of $\mathbb{S}^m$ is 
$Vol_{S^m}=\sum_{i} (-1)^{i-1}\phi_i d\phi^{1\ldots\hat{i}\ldots m} $, 
and $*\xi$ is the unique 2-form s.t. $\xi\wedge * \xi= \|\xi\|^2Vol_{S^m}$.
Using (7) we see that $\|\xi\|=\|*\xi\|=1$.
Hence
$$\begin{array}{l}
\xi= \phi_1d\phi^2- \phi_2d\phi^1+ \phi_3d\phi^4- \phi_4d\phi^3\\
*\xi= d\phi^1\wedge d\phi^2+ d\phi^3\wedge d\phi^4=\frac{1}{2}
d\xi=:
d*\omega. \end{array}$$
Therefore, we may take
$*\omega= \frac{1}{2}\xi$, that is
$$\omega=\frac{1}{2}*\xi=\frac{1}{2}(d\phi^1\wedge d\phi^2 +
d\phi^3\wedge d\phi^4)=
\frac{1}{2}\phi^*\varpi. $$
Hence, to prove Theorem 1.2  and Corollary 1.1
we have to verify that, for any functions 
$f,h\in C^{\infty}
(\mathbb{S}^3)$, one of the following equivalent inequalities holds:
\begin{eqnarray}
\int_{\mathbb{S}^3} -3\omega(\nabla f,\nabla h)dM
=\int_{\mathbb{S}^3} -3f\xi(\nabla h)dM\leq \|\nabla f\|_{L^2}
\|\nabla h\|_{L^2}\\
\int_{\mathbb{S}^3} -6\omega(\nabla f,\nabla h)dM
=\int_{\mathbb{S}^3} -6f\xi(\nabla h)dM\leq \|\nabla f\|_{L^2}^2+
\|\nabla h\|_{L^2}^2.\nonumber
\end{eqnarray}
By Theorem 1.1 we  only need to consider both
$f,h\in E_{\lambda_l}$, for some $l$.
Note that  $\lambda_3=15$ and  since $\Omega$ is a calibration, 
$ \|\xi(X)\|\leq \|X\|$.\\[2mm]
{\bf Lemma 4.1.} \em 
If  $f, h\in E_{\lambda_3}^+$ are nonzero, 
(8) holds, with strict inequality. \em \\[2mm]
{\bf  Proof.}
By Schwartz inequality and Rayleigh characterization
$$\int_{\mathbb{S}^3} -3f\xi(\nabla h)dM\leq 3\|f\|_{L^2}
\|\nabla h\|_{L^2}\leq \frac{3}{\sqrt{\lambda_3}}
\|\nabla f\|_{L^2}\|\nabla h\|_{L^2} <
\|\nabla f\|_{L^2}\|\nabla h\|_{L^2}, $$
with strict inequality in the last one, since neither $f$ nor $h$ 
may be constant.\qed\par

We now verify that (8) holds for $f,h \in E_{\lambda_1}$ and 
$f,h\in E_{\lambda_2}$.
From (7) and Lemma 3.1, we have for $i\neq j$
\begin{eqnarray} 
\begin{array}{lcl}
\int_{\mathbb{S}^3}\phi^2dM=\frac{1}{4}|\mathbb{S}^3|,&& 
\int_{\mathbb{S}^3}\phi^2_i\phi_j^2dM=
\frac{1}{ 6}\int_{\mathbb{S}^3}\phi^2dM \\
\int_{\mathbb{S}^3}\phi^4dM=\frac{1}{2}\int_{\mathbb{S}^3}\phi^2dM,&&
\int_{\mathbb{S}^3}\|\nabla \phi\|^2dM= 3\int_{\mathbb{S}^3}\phi^2dM\\[1mm]
\omega(\nabla \phi_1,\nabla \phi_2) =\frac{1}{2}(1-\phi_1^2-\phi_2^2) 
&&
\omega(\nabla \phi_1,\nabla \phi_3) =\frac{1}{2}(-\phi_2\phi_3+\phi_1\phi_4)
\\ 
\omega(\nabla \phi_1,\nabla \phi_4) =\frac{1}{2}(-\phi_2\phi_4-\phi_1\phi_3)
&&
\omega(\nabla \phi_2,\nabla \phi_3) =\frac{1}{2}(\phi_1\phi_3+\phi_4\phi_2)
\\
\omega(\nabla \phi_2,\nabla \phi_4) =\frac{1}{2}(\phi_1\phi_4-\phi_2\phi_3)
& &\omega(\nabla \phi_3,\nabla \phi_4) =\frac{1}{2}(1-\phi_3^2-\phi_4^2).
\end{array}
\end{eqnarray}
and moreover\\[2mm]
{\bf Lemma 4.2.} 
$$\begin{array}{l}
3\int\omega(\nabla\phi_1,\nabla\phi_2)=3\int\phi^2=\|\nabla \phi_1\|_{L^2}
\|\nabla \phi_2\|_{L^2}=\|\nabla \phi\|_{L^2}^2\\
3\int\omega(\nabla\phi_3,\nabla\phi_4)=3\int\phi^2
=\|\nabla \phi_3\|_{L^2}
\|\nabla \phi_4\|_{L^2}=\|\nabla \phi\|_{L^2}^2\\
-3\int\omega(\nabla\phi_i,\nabla\phi_j)=0~~~\mbox{for other }ij\\[1mm]
-3\int\phi_k\omega(\nabla\phi_i,\nabla\phi_j)=0~~~\forall i,j,k\\
-3\int\phi_1^2\omega(\nabla\phi_1,\nabla\phi_2)=
-3\int\phi_2^2\omega(\nabla\phi_1,\nabla\phi_2)=
-\frac{1}{2}\int \phi^2 \\
-3\int\phi_3^2\omega(\nabla\phi_1,\nabla\phi_2)=
-3\int\phi_4^2\omega(\nabla\phi_1,\nabla\phi_2)=
-\int \phi^2 \\
-3\int\phi_1^2\omega(\nabla\phi_3,\nabla\phi_4)=
-3\int\phi_2^2\omega(\nabla\phi_3,\nabla\phi_4)=
-\int \phi^2 \\
-3\int\phi_3^2\omega(\nabla\phi_3,\nabla\phi_4)=
-3\int\phi_4^2\omega(\nabla\phi_3,\nabla\phi_4)=
-\frac{1}{2}\int \phi^2 \\[2mm]
-3\int\phi_1\phi_4\omega(\nabla\phi_1,\nabla\phi_3)=
-3\int\phi_1\phi_3\omega(\nabla\phi_2,\nabla\phi_3)=
-\frac{1}{4}\int\phi^2\\
-3\int\phi_1\phi_3\omega(\nabla\phi_1,\nabla\phi_4)=
-3\int\phi_2\phi_3\omega(\nabla\phi_2,\nabla\phi_4)
=\frac{1}{4}\int\phi^2\\
-3\int\phi_2\phi_3\omega(\nabla\phi_1,\nabla\phi_3)=
-3\int\phi_2\phi_4\omega(\nabla\phi_1,\nabla\phi_4)=
\frac{1}{4}\int\phi^2\\
-3\int\phi_2\phi_4\omega(\nabla\phi_2,\nabla\phi_3)=
-3\int\phi_1\phi_4\omega(\nabla\phi_2,\nabla\phi_4)=
-\frac{1}{4}\int\phi^2\\
-3\int\phi_i\phi_j\omega(\nabla\phi_k,\nabla\phi_s)=0~~\mbox{for other cases.}
\\
\end{array}
$$
\noindent
{\bf Lemma 4.3.} \em 
If  $f, h\in E_{\lambda_1}$,  that is $f=\sum_i\mu_i\phi_i$, 
$h=\sum_j\sigma_j\phi_j$,  for some constant $\mu_i, \sigma_j$,
then (8) holds,
with equality if and only if $\sigma_2=-\mu_1$, $\sigma_1=\mu_2$, 
$\sigma_4=-\mu_3$, $\sigma_3= \mu_4$.
\em \\[2mm]
{\bf  Proof. }
Using the previous lemma, 
$$\begin{array}{lcl}
-3\int\omega(\nabla f,\nabla h)dM
&=&
(\mu_1\sigma_2-\mu_2\sigma_1)\int -3\omega(\nabla \phi_1,\nabla \phi_2)
+(\mu_3\sigma_4-\mu_4\sigma_3)\int -3\omega(\nabla \phi_3,\nabla \phi_4)\\
&=& -(\mu_1\sigma_2-\mu_2\sigma_1+\mu_3\sigma_4-\mu_4\sigma_3)
\|\nabla \phi\|_{L^2}^2\\
&\leq& \frac{1}{2}(\sum_i \mu_i^2+ \sigma_i^2)\|\nabla \phi\|_{L^2}^2
=\frac{1}{2}(\|\nabla f\|_{L^2}^2+ \|\nabla h\|_{L^2}^2). 
\end{array} $$
The equality case follows immediately. \qed \\[2mm]
{\bf Lemma 4.4.} \em 
If  $f, h\in E_{\lambda_2}$  are nonzero,
then (8) holds
with strict inequality. \em \\[2mm]
{\bf Proof. } Set
$f=\sum_i\alpha_i\phi_i^2+ \sum_{i<j}A_{ij}\phi_i\phi_j$,
 and $h=\sum_i\beta_i\phi_i^2+ \sum_{i<j}B_{ij}\phi_i\phi_j$, where 
$\alpha_i, A_{ij}$, $\beta_i, B_{ij}$ are constants.
Now we compute
{\small \begin{eqnarray*}
\lefteqn{-3\int \omega(\nabla f, \nabla h)=}\\
&&-3\int \omega(\nabla \phi_1, \nabla \phi_2)
[(2\alpha_1\phi_1+ A_{12}\phi_2+A_{13}\phi_3+A_{14}\phi_4)
(2\beta_2\phi_2+ B_{12}\phi_1+B_{23}\phi_3+B_{24}\phi_4)\\
&&\quad \quad \quad -(2\alpha_2\phi_2+ A_{12}\phi_1+A_{23}\phi_3+A_{24}\phi_4)
(2\beta_1\phi_1+ B_{12}\phi_2+B_{13}\phi_3+B_{14}\phi_4)]\\
&&-3\int \omega(\nabla \phi_1, \nabla \phi_3)
[(2\alpha_1\phi_1+ A_{12}\phi_2+A_{13}\phi_3+A_{14}\phi_4)
(2\beta_3\phi_3+ B_{13}\phi_1+B_{23}\phi_2+B_{34}\phi_4)\\
&&\quad \quad \quad -(2\alpha_3\phi_3+ A_{13}\phi_1+A_{23}\phi_2+A_{34}\phi_4)
(2\beta_1\phi_1+ B_{12}\phi_2+B_{13}\phi_3+B_{14}\phi_4)]\\
&&-3\int \omega(\nabla \phi_1, \nabla \phi_4)
[(2\alpha_1\phi_1+ A_{12}\phi_2+A_{13}\phi_3+A_{14}\phi_4)
(2\beta_4\phi_4+ B_{14}\phi_1+B_{24}\phi_2+B_{34}\phi_3)\\
&&\quad \quad \quad -(2\alpha_4\phi_4+ A_{14}\phi_1+A_{24}\phi_2+A_{34}\phi_3)
(2\beta_1\phi_1+ B_{12}\phi_2+B_{13}\phi_3+B_{14}\phi_4)]\\
&&-3\int \omega(\nabla \phi_2, \nabla \phi_3)
[(2\alpha_2\phi_2+ A_{12}\phi_1+A_{23}\phi_3+A_{24}\phi_4)
(2\beta_3\phi_3+ B_{13}\phi_1+B_{23}\phi_2+B_{34}\phi_4)\\
&&\quad \quad \quad -(2\alpha_3\phi_3+ A_{13}\phi_1+A_{23}\phi_2+A_{34}\phi_4)
(2\beta_2\phi_2+ B_{12}\phi_1+B_{24}\phi_4+B_{23}\phi_3)]\\
&&-3\int \omega(\nabla \phi_2, \nabla \phi_4)
[(2\alpha_2\phi_2+ A_{12}\phi_1+A_{23}\phi_3+A_{24}\phi_4)
(2\beta_4\phi_4+ B_{14}\phi_1+B_{24}\phi_2+B_{34}\phi_3)\\
&&\quad \quad \quad -(2\alpha_4\phi_4+ A_{14}\phi_1+A_{24}\phi_2+A_{34}\phi_3)
(2\beta_2\phi_2+ B_{12}\phi_1+B_{24}\phi_4+B_{23}\phi_3)]\\
&&-3\int \omega(\nabla \phi_3, \nabla \phi_4)
[(2\alpha_3\phi_3+ A_{13}\phi_1+A_{23}\phi_2+A_{34}\phi_4)
(2\beta_4\phi_4+ B_{14}\phi_1+B_{24}\phi_2+B_{34}\phi_3)\\
&&\quad \quad \quad -(2\alpha_4\phi_4+ A_{14}\phi_1+A_{24}\phi_2+A_{34}\phi_3)
(2\beta_3\phi_3+ B_{13}\phi_1+B_{23}\phi_2+B_{34}\phi_4)].
\end{eqnarray*}}
Thus, using Lemma 4.2, 
{\small \begin{eqnarray*}
\lefteqn{-3\int \omega(\nabla f, \nabla h)=}\\
&-3\int \omega(\nabla \phi_1, \nabla \phi_2)&
[2\alpha_1B_{12}\phi_1^2+ 2\beta_2A_{12}\phi_2^2
+A_{13}B_{23}\phi_3^2+A_{14}B_{24}\phi_4^2\\
&&-2\beta_1A_{12}\phi_1^2 -2\alpha_2B_{12}\phi_2^2
-A_{23}B_{13}\phi_3^2-A_{24}B_{14}\phi_4^2]\\
&-3\int \omega(\nabla \phi_3, \nabla \phi_4)&
[A_{13}B_{14}\phi_1^2+ A_{23}B_{24}\phi_2^2
+2\alpha_{3}B_{34}\phi_3^2+2\beta_{4}A_{34}\phi_4^2\\
&&-A_{14}B_{13}\phi_1^2 -A_{24}B_{23}\phi_2^2
-2\beta_{3}A_{34}\phi_3^2-2\alpha_{4}B_{34}\phi_4^2]\\
&-3\int \omega(\nabla \phi_1, \nabla \phi_3)&
[2\alpha_1B_{34}\phi_1\phi_4+ A_{14}B_{13}\phi_1\phi_4
-A_{13}B_{14}\phi_1\phi_4-2\beta_1A_{34}\phi_1\phi_4\\
&&+2\beta_3A_{12}\phi_2\phi_3 +A_{13}B_{23}\phi_2\phi_3
-A_{23}B_{13}\phi_2\phi_3-2\alpha_3B_{12}\phi_2\phi_3]\\
&-3\int \omega(\nabla \phi_1, \nabla \phi_4)&
[2\alpha_1B_{34}\phi_1\phi_3+ A_{13}B_{14}\phi_1\phi_3
-A_{14}B_{13}\phi_1\phi_3-2\beta_1A_{34}\phi_1\phi_3\\
&&+2\beta_4A_{12}\phi_2\phi_4 +A_{14}B_{24}\phi_2\phi_4
-A_{24}B_{14}\phi_2\phi_4-2\alpha_4B_{12}\phi_2\phi_4]\\
&-3\int \omega(\nabla \phi_2, \nabla \phi_3)&
[2\beta_3A_{12}\phi_1\phi_3+ A_{23}B_{13}\phi_1\phi_3
-A_{13}B_{23}\phi_1\phi_3-2\alpha_3B_{12}\phi_1\phi_3\\
&&+2\alpha_2B_{34}\phi_2\phi_4 +A_{24}B_{23}\phi_2\phi_4
-A_{23}B_{24}\phi_2\phi_4-2\beta_2A_{34}\phi_2\phi_4]\\
&-3\int \omega(\nabla \phi_2, \nabla \phi_4)&
[2\beta_4A_{12}\phi_1\phi_4+ A_{24}B_{14}\phi_1\phi_4
-A_{14}B_{24}\phi_1\phi_4-2\alpha_4B_{12}\phi_1\phi_4\\
&&+2\alpha_2B_{34}\phi_2\phi_3 +A_{23}B_{24}\phi_2\phi_3
-A_{24}B_{23}\phi_2\phi_3-2\beta_2A_{34}\phi_2\phi_3]
\end{eqnarray*}}
{\small 
\begin{eqnarray*}
\begin{array}{rll}
=&\int \phi^2\La{\{}&
-\frac{1}{2}[2\alpha_1B_{12}+2\beta_2A_{12}-2\beta_1A_{12}-2\alpha_2B_{12}
+2\alpha_3B_{34}+2\beta_4A_{34}-2\beta_3A_{34}-2\alpha_4B_{34}]
\\
&&-[A_{13}B_{23}+A_{14}B_{24}-A_{23}B_{13}-A_{24}B_{14}+
A_{13}B_{14}+A_{23}B_{24}-A_{14}B_{13}-A_{24}B_{23}]\\
&&+\frac{1}{4}\La{[}-2\alpha_1B_{34}-A_{14}B_{13}+A_{13}B_{14}+ 2\beta_1A_{34}
+ 2\beta_3A_{12}+A_{13}B_{23}-A_{23}B_{13}-2\alpha_3B_{12} \\[1mm]
&&\quad\quad+2\alpha_1B_{34}+A_{13}B_{14}-A_{14}B_{13}- 2\beta_1A_{34}
+ 2\beta_4A_{12}+A_{14}B_{24}-A_{24}B_{14}-2\alpha_4B_{12}\\[1mm]
&&\quad\quad-2\beta_3A_{12}-A_{23}B_{13}+A_{13}B_{23}+ 2\alpha_3B_{12}
- 2\alpha_2B_{34}-A_{24}B_{23}+A_{23}B_{24}+2\beta_2A_{34}\\[1mm]
&&\quad\quad -2\beta_41A_{12}-A_{24}B_{14}+A_{14}B_{24}+ 2\alpha_4B_{12}
+ 2\alpha_2B_{34}+A_{23}B_{24}-A_{24}B_{23}-2\beta_2A_{34}\La{]}~~\La{\}} 
\\[2mm]
=&\int\phi^2\La{\{}&-[\alpha_1B_{12}+\beta_2A_{12}-\beta_1A_{12}-\alpha_2B_{12}
+\alpha_3B_{34}+\beta_4A_{34}-\beta_3A_{34}-\alpha_4B_{34}] \\
&&-[A_{13}B_{23}+A_{14}B_{24}-A_{23}B_{13}-A_{24}B_{14}+
A_{13}B_{14}+A_{23}B_{24}-A_{14}B_{13}-A_{24}B_{23}]\\
&&+\frac{1}{2}\La{[}
-A_{14}B_{13}+A_{13}B_{14}+A_{13}B_{23}-A_{23}B_{13}
+A_{14}B_{24}-A_{24}B_{14}-A_{24}B_{23}+A_{23}B_{24}\La{]}~~\La{\}}\\[2mm]
=&\int\phi^2\La{\{}&
[-\alpha_1B_{12}-\beta_2A_{12}+\beta_1A_{12}+\alpha_2B_{12}
-\alpha_3B_{34}-\beta_4A_{34}+\beta_3A_{34}+\alpha_4B_{34}]\\
&&+ \frac{1}{2}[-A_{13}B_{23}-A_{14}B_{24}+A_{23}B_{13}+A_{24}B_{14}-
A_{13}B_{14}-A_{23}B_{24}+A_{14}B_{13}+A_{24}B_{23}]~~ \La{\}}
\end{array}
\end{eqnarray*}}
and applying the same lemmas we see that
$$\|\nabla f\|_{L^2}^2=\left
[2(\sum_k \alpha_k^2)-\frac{4}{3}(\sum_{i<j} \alpha_i
\alpha_j)+ \frac{4}{3}(\sum_{i<j} A^2_{ij})\right]\int \phi^2.$$
Hence, we have to verify if  the following inequality is true:
\begin{eqnarray}
[-\alpha_1B_{12}-\beta_2A_{12}+\beta_1A_{12}+\alpha_2B_{12}
-\alpha_3B_{34}-\beta_4A_{34}+\beta_3A_{34}+\alpha_4B_{34}]\label{10}\\
+\frac{1}{2} [-A_{13}B_{23}-A_{14}B_{24}+A_{23}B_{13}+A_{24}B_{14}-
A_{13}B_{14}-A_{23}B_{24}+A_{14}B_{13}+A_{24}B_{23}]\label{11}\\
+\frac{2}{3}(\sum_{i<j}\alpha_i\alpha_j + \beta_i\beta_j)\label{12}\\
\leq \sum_k (\alpha_k^2 + \beta_k^2) 
+ \frac{2}{3}(\sum_{i<j} A^2_{ij}+ B_{ij}^2).\label{13}
\end{eqnarray}
This is equivalent to prove the inequalities
\begin{eqnarray}
(\ref{11})&\leq & \frac{2}{3}
(A_{13}^2+A_{14}^2+A_{23}^2+A_{24}^2+B_{13}^2+B_{14}^2+B_{23}^2+B_{24}^2)
\label{14}\\[2mm]
(\ref{10})+(\ref{12})&\leq &  \sum_k (\alpha_k^2 + \beta_k^2) + \frac{2}{3}(
A_{12}^2+A_{34}^2+B_{12}^2+B_{34}^2).\label{15}
\end{eqnarray}
Note that
\begin{eqnarray*}
2\times(\ref{11})
&\leq & (A_{13}^2+A_{14}^2+A_{23}^2+A_{24}^2+B_{13}^2
+B_{14}^2+B_{23}^2+B_{24}^2)\\
&\leq& \frac{4}{3}(A_{13}^2+A_{14}^2+A_{23}^2+A_{24}^2+B_{13}^2
+B_{14}^2+B_{23}^2+B_{24}^2),
\end{eqnarray*}
and so inequality (\ref{14}) holds, with equality if and only if 
$$A_{13}=A_{14}=A_{23}=A_{24}=B_{13}=B_{14}=B_{23}=B_{24}=0.$$
Now
\begin{eqnarray}
3\times (\ref{10})
&=& 3(\alpha_2-\alpha_1)B_{12}- 3(\beta_2-\beta_1)A_{12}
+ 3(\alpha_4-\alpha_3)B_{34}+ 3(-\beta_4+\beta_3)A_{34} \nonumber \\
&\leq& 
 \frac{3}{2}( (\alpha_2-\alpha_1)^2+ (\beta_2-\beta_1)^2
+ (\alpha_4-\alpha_3)^2+ (-\beta_4+\beta_3)^2)\nonumber \\
&&+  \frac{3}{2}(A_{12}^2+ A_{34}^2+B_{12}^2+B_{34}^2) \nonumber \\
&\leq &\frac{3}{2}( (\alpha_2-\alpha_1)^2+ (\beta_2-\beta_1)^2
+ (\alpha_4-\alpha_3)^2+ (-\beta_4+\beta_3)^2)\label{16}\\
&&+ 2(A_{12}^2+ A_{34}^2+B_{12}^2+B_{34}^2).\label{17}
\end{eqnarray}
We will prove that 
\begin{eqnarray} \label{18}
(\ref{16})+  3\times(\ref{12})
&\leq & \sum_k3(\alpha_k^2+\beta_k^2),
\end{eqnarray}
with equality iff $\alpha_1=\alpha_2=\alpha_3=\alpha_4$ and
$\beta_1=\beta_2=\beta_3=\beta_4$,
which proves that  (\ref{15}) holds. Furthermore, from (\ref{17})  we see
that equality in (\ref{15})  is achieved iff
$$A_{12}= A_{34}=B_{12}=B_{34}=0,\mbox{~~~
and for all~} i,j ~~~\alpha_i=\alpha_j,~~~\beta_i=\beta_j.$$
In order to prove  (\ref{18}) we only have to show that
$$\frac{3}{2}((\alpha_2-\alpha_1)^2+ (\alpha_4-\alpha_3)^2)
+2\sum_{i<j}\alpha_i\alpha_j\leq 3\sum_k\alpha_k^2,$$
or equivalently, that
$$-2\alpha_1\alpha_2-2\alpha_3\alpha_4+4\alpha_1\alpha_3+4\alpha_1\alpha_4
+4\alpha_2\alpha_3+4\alpha_2\alpha_4\leq 3\sum_k \alpha_k^2.$$
But this is just
$$(\alpha_1-\alpha_3)^2+(\alpha_3-\alpha_2)^2+ (\alpha_2-\alpha_4)^2
+(\alpha_4-\alpha_1)^2+(\alpha_1+\alpha_2-\alpha_3-\alpha_4)^2\geq 0,$$
with equality to zero iff $\alpha_i=\alpha_j$ $\forall ij$. 
We have proved that inequality (8) is satisfied,
 with equality iff
$f=\alpha(\sum_k\phi_k^2)=  \alpha$ constant and $ h$ constant, and so
they must vanish. \qed\\[1mm]

Theorem 1.1,  with Lemmas 4.1, 4.3
 and 4.4, prove that
(8) holds for any pair of functions $(f,h)$, and so   
Theorem 1.2 is proved.
Corollary 1.1 follows from these lemmas.\par

In  (Salavessa, 2010, Theorem 4.2) a
 uniqueness theorem was obtained, on a class of closed
$m$-dimensional submanifolds with parallel mean curvature and 
calibrated extended tangent
in a Euclidean space $\mathbb{R}^{m+n}$,  and satisfying an integral height inequality.
We will recall such results for the case $\Omega$ parallel.
We denote by $B^{\nu}$  the $\nu$-component of the second fundamental form $B$
and by $B^F$ the $F$-component, $B= B^{\nu}+ B^F$, where $F$ is the
orthogonal complement of $\nu$ in the normal bundle.\\[2mm]
{\bf Theorem 4.1.} \em  If $\Omega$ is a parallel calibration of rank $(m+1)$
on
$\mathbb{R}^{m+n}$,  and $\phi:M\to \mathbb{R}^{m+n}$ is an immersed
 closed $\Omega$-stable $m$-dimensional submanifold with parallel 
mean curvature and calibrated extended tangent space, and
\begin{equation}
\int_M S(2 + h \|H\|)dM \leq 0, 
\end{equation}
where 
$ h=\langle \phi, \nu\rangle$ and
$S= \sum_{ij} \langle \phi, (B(e_i,e_j))^F \rangle B^{\nu}(e_i,e_j),$
then $\phi$ is pseudo-umbilical and $S=0$. Furthermore, if $NM$ is
a trivial bundle, then the minimal calibrated extension of $M$ is a
Euclidean space $\mathbb{R}^{m+1}$, and $M$ is a Euclidean $m$-sphere. \em 
\\[1mm]
\noindent 
 Theorem 1.3 is an immediate consequence of 
Theorem 1.2 and the above theorem.\\ \par
\textbf{Acknowledgements} \\[-2mm]

The author would like to thank
Dr.\ Ana Cristina
Ferreira  and the Universidade do Minho, Braga, for their
   hospitality  during the \em Third
Minho Meeting on Mathematical Physics \em  at Centro de Matem\'{a}tica da UM, 
in November 2011,
where the final part of this work was completed. \\
\par
\textbf{References}\\[3mm]
{\small J.L. Barbosa and M. do Carmo (1984). {\sl Stability of 
hypersurfaces with constant mean curvature},  Math.\ Z., 185(3), 339-353.
\\[1mm]
{\small J.L. Barbosa, M. do Carmo and J. Eschenburg (1988). {\sl 
Stability of hypersurfaces of constant mean curvature in Riemannian manifolds},
 Math.\ Z., 197(1), 123- 138.
\\[1mm]
{\small F. Morgan (2000), {\sl Perimeter-minimizing curves and surfaces in 
$\mathbb{R}^n$  enclosing prescribed multi-volume}, Asian J.\ Math., 4 (2), 
 373-382.}\\[1mm]
{\small I.M.C.\ Salavessa (2010). {\sl Stability of submanifolds with parallel
mean curvature in calibrated manifolds},  
Bull.\  Braz.\ Math.\ Soc.\, NS,  41(4),  495-530. }
\\[1mm] 
{\small I.M.C.\ Salavessa \& A.\ Pereira do Vale (2006). {\sl
Transgression forms in dimension 4}, IJGMMP, 3(4-5), 1221-1254.}
\end{document}